\documentclass{amsart}
\usepackage{amsfonts}

\newtheorem{theorem}{Theorem}
\theoremstyle{plain}

\newtheorem{corollary}{Corollary}

\newtheorem{lemma}{Lemma}

\newtheorem{proposition}{Proposition}
\newtheorem{remark}{Remark}

\numberwithin{equation}{section}
\input{tcilatex}

\begin{document}
\title[Schwarz Related Inequalities]{A Potpourri of Schwarz Related
Inequalities in Inner Product Spaces}
\author{S.S. Dragomir}
\address{School of Computer Science and Mathematics\\
Victoria University of Technology\\
PO Box 14428, MCMC 8001\\
VIC, Australia.}
\email{sever@matilda.vu.edu.au}
\urladdr{http://rgmia.vu.edu.au/SSDragomirWeb.html}
\date{30 September, 2004.}

\begin{abstract}
In this paper we obtain some new Schwarz related inequalities in inner
product spaces over the real or complex number field. Applications for the
generalized triangle inequality are also given.
\end{abstract}

\maketitle

\section{Introduction}

Let $\left( H;\left\langle \cdot ,\cdot \right\rangle \right) $ be an inner
product space over the real or complex number field $\mathbb{K}$. One of the
most important inequalities in inner product spaces with numerous
applications, is the Schwarz inequality; that may be written in two forms: 
\begin{equation}
\left\vert \left\langle x,y\right\rangle \right\vert ^{2}\leq \left\Vert
x\right\Vert ^{2}\left\Vert y\right\Vert ^{2},\ \ \ \ x,y\in H\qquad \text{%
(quadratic form)}  \label{1.1}
\end{equation}%
or, equivalently, 
\begin{equation}
\left\vert \left\langle x,y\right\rangle \right\vert \leq \left\Vert
x\right\Vert \left\Vert y\right\Vert ,\ \ \ \ x,y\in H\qquad \text{(simple
form).}  \label{1.1b}
\end{equation}%
The case of equality holds in (\ref{1.1}) (or (\ref{1.1b})) if and only if
the vectors $x$ and $y$ are linearly dependent.

In 1966, S. Kurepa \cite{K}, gave the following refinement of the quadratic
form for the complexification of a real inner product space:

\begin{theorem}
\label{t1.1}Let $\left( H;\left\langle \cdot ,\cdot \right\rangle \right) $
be a real Hilbert space and $\left( H_{\mathbb{C}},\left\langle \cdot ,\cdot
\right\rangle _{\mathbb{C}}\right) $ the complexification of $H.$ Then for
any pair of vectors $a\in H,$ $z\in H_{\mathbb{C}}$%
\begin{equation}
\left\vert \left\langle z,a\right\rangle _{\mathbb{C}}\right\vert ^{2}\leq 
\frac{1}{2}\left\Vert a\right\Vert ^{2}\left( \left\Vert z\right\Vert _{%
\mathbb{C}}^{2}+\left\vert \left\langle z,\bar{z}\right\rangle _{\mathbb{C}%
}\right\vert \right) \leq \left\Vert a\right\Vert ^{2}\left\Vert
z\right\Vert _{\mathbb{C}}^{2}.  \label{1.2}
\end{equation}
\end{theorem}

In 1985, S.S. Dragomir \cite[Theorem 2]{SSD1} obtained a refinement of the
simple form of Schwarz inequality as follows:

\begin{theorem}
\label{t1.2}Let $\left( H;\left\langle \cdot ,\cdot \right\rangle \right) $
be a real or complex inner product space and $x,y,e\in H$ with $\left\Vert
e\right\Vert =1.$ Then we have the inequality 
\begin{equation}
\left\Vert x\right\Vert \left\Vert y\right\Vert \geq \left\vert \left\langle
x,y\right\rangle -\left\langle x,e\right\rangle \left\langle
e,y\right\rangle \right\vert +\left\vert \left\langle x,e\right\rangle
\left\langle e,y\right\rangle \right\vert \geq \left\vert \left\langle
x,y\right\rangle \right\vert .  \label{1.3}
\end{equation}
\end{theorem}

For other similar results, see \cite{DS} and \cite{DM}.

A refinement of the \textit{weaker version }of the Schwarz inequality, i.e.,%
\begin{equation}
\func{Re}\left\langle x,y\right\rangle \leq \left\Vert x\right\Vert
\left\Vert y\right\Vert ,\ \ \ \ x,y\in H  \label{1.5}
\end{equation}%
has been established in \cite{SSD2}:

\begin{theorem}
\label{t1.3}Let $\left( H;\left\langle \cdot ,\cdot \right\rangle \right) $
be a real or complex inner product space and $x,y,e\in H$ with $\left\Vert
e\right\Vert =1.$ \ If $r_{1},r_{2}>0$ and $x,y\in H$ are such that%
\begin{equation}
\left\Vert x-y\right\Vert \geq r_{2}\geq r_{1}\geq \left\vert \left\Vert
x\right\Vert -\left\Vert y\right\Vert \right\vert   \label{1.6}
\end{equation}%
then we have the following refinement of the weak Schwarz inequality%
\begin{equation}
\left\Vert x\right\Vert \left\Vert y\right\Vert -\func{Re}\left\langle
x,y\right\rangle \geq \frac{1}{2}\left( r_{2}^{2}-r_{1}^{2}\right) \qquad
\left( \geq 0\right) .  \label{1.7}
\end{equation}%
The constant $\frac{1}{2}$ is best possible in the sense that it cannot be
replaced by a larger quantity.
\end{theorem}

For other recent results see the paper mentioned above, \cite{SSD2}.

In practice, one may need reverses of the Schwarz inequality, namely, upper
bounds for the quantities%
\begin{equation*}
\left\Vert x\right\Vert \left\Vert y\right\Vert -\func{Re}\left\langle
x,y\right\rangle ,\qquad \left\Vert x\right\Vert ^{2}\left\Vert y\right\Vert
^{2}-\left( \func{Re}\left\langle x,y\right\rangle \right) ^{2}
\end{equation*}%
and%
\begin{equation*}
\frac{\left\Vert x\right\Vert \left\Vert y\right\Vert }{\func{Re}%
\left\langle x,y\right\rangle }
\end{equation*}%
or the coresponding expressions where $\func{Re}\left\langle
x,y\right\rangle $ is repalced by either $\left\vert \func{Re}\left\langle
x,y\right\rangle \right\vert $ or $\left\vert \left\langle x,y\right\rangle
\right\vert ,$ under suitable assumptions for the vectors $x,y$ in an inner
product space $\left( H;\left\langle \cdot ,\cdot \right\rangle \right) $
over the real or complex number field $\mathbb{K}$.

In this class of results, we mention the following recent reverses of the
Schwarz inequality due to the present author, that can be found, for
instance, in the survey work \cite{SSD3}, where more specific references are
provided:

\begin{theorem}
\label{t1.4}Let $\left( H;\left\langle \cdot ,\cdot \right\rangle \right) $
be an inner product space over $\mathbb{K}$ $\left( \mathbb{K}=\mathbb{C},%
\mathbb{R}\right) .$ If $a,A\in \mathbb{K}$ and $x,y\in H$ are such that
either%
\begin{equation}
\func{Re}\left\langle Ay-x,x-ay\right\rangle \geq 0,  \label{1.8}
\end{equation}%
or, equivalently,%
\begin{equation}
\left\Vert x-\frac{A+a}{2}y\right\Vert \leq \frac{1}{2}\left\vert
A-a\right\vert \left\Vert y\right\Vert ,  \label{1.9}
\end{equation}%
then the following reverse for the quadratic form of the Schwarz inequality%
\begin{align}
(0& \leq )\left\Vert x\right\Vert ^{2}\left\Vert y\right\Vert
^{2}-\left\vert \left\langle x,y\right\rangle \right\vert ^{2}  \label{1.10}
\\
& \leq \left\{ 
\begin{array}{l}
\frac{1}{4}\left\vert A-a\right\vert ^{2}\left\Vert y\right\Vert
^{4}-\left\vert \frac{A+a}{2}\left\Vert y\right\Vert ^{2}-\left\langle
x,y\right\rangle \right\vert ^{2} \\ 
\\ 
\frac{1}{4}\left\vert A-a\right\vert ^{2}\left\Vert y\right\Vert
^{4}-\left\Vert y\right\Vert ^{2}\func{Re}\left\langle Ay-x,x-ay\right\rangle%
\end{array}%
\right.  \notag \\
& \leq \frac{1}{4}\left\vert A-a\right\vert ^{2}\left\Vert y\right\Vert ^{4}
\notag
\end{align}%
holds.

If in addition, we have $\func{Re}\left( A\bar{a}\right) >0,$ then%
\begin{equation}
\left\Vert x\right\Vert \left\Vert y\right\Vert \leq \frac{1}{2}\cdot \frac{%
\func{Re}\left[ \left( \bar{A}+\bar{a}\right) \left\langle x,y\right\rangle %
\right] }{\sqrt{\func{Re}\left( A\bar{a}\right) }}\leq \frac{1}{2}\cdot 
\frac{\left\vert A+a\right\vert }{\sqrt{\func{Re}\left( A\bar{a}\right) }}%
\left\vert \left\langle x,y\right\rangle \right\vert ,  \label{1.11}
\end{equation}%
and%
\begin{equation}
(0\leq )\left\Vert x\right\Vert ^{2}\left\Vert y\right\Vert ^{2}-\left\vert
\left\langle x,y\right\rangle \right\vert ^{2}\leq \frac{1}{4}\cdot \frac{%
\left\vert A-a\right\vert ^{2}}{\func{Re}\left( A\bar{a}\right) }\left\vert
\left\langle x,y\right\rangle \right\vert ^{2}.  \label{1.12}
\end{equation}%
Also, if (\ref{1.8}) or (\ref{1.9}) are valid and $A\neq -a,$ then we have
the reverse for the simple form of Schwarz inequality%
\begin{align}
(0& \leq )\left\Vert x\right\Vert \left\Vert y\right\Vert -\left\vert
\left\langle x,y\right\rangle \right\vert \leq \left\Vert x\right\Vert
\left\Vert y\right\Vert -\left\vert \func{Re}\left[ \frac{\bar{A}+\bar{a}}{%
\left\vert A+a\right\vert }\left\langle x,y\right\rangle \right] \right\vert
\label{1.13} \\
& \leq \left\Vert x\right\Vert \left\Vert y\right\Vert -\func{Re}\left[ 
\frac{\bar{A}+\bar{a}}{\left\vert A+a\right\vert }\left\langle
x,y\right\rangle \right] \leq \frac{1}{4}\cdot \frac{\left\vert
A-a\right\vert ^{2}}{\left\vert A+a\right\vert }\left\Vert y\right\Vert ^{2}.
\notag
\end{align}%
The multiplicative constants $\frac{1}{4}$ and $\frac{1}{2}$ above are best
possible in the sense that they cannot be replaced by a smaller quantity.
\end{theorem}

For some classical results related to Schwarz inequality, see \cite{B}, \cite%
{FK}, \cite{PE}, \cite{P}, \cite{R} and the references therein.

The main aim of the present paper is to point out other results in
connection with both the quadratic and simple forms of the Schwarz
inequality. As applications, some reverse results for the generalised
triangle inequality, i.e., upper bounds for the quantity%
\begin{equation*}
(0\leq )\sum_{i=1}^{n}\left\Vert x_{i}\right\Vert -\left\Vert
\sum_{i=1}^{n}x_{i}\right\Vert
\end{equation*}%
under various assumptions for the vectors $x_{i}\in H,$ $i\in \left\{
1,\dots ,n\right\} ,$ are established.

\section{Inequalities Related to Schwarz's One\label{s2}}

The following result holds.

\begin{proposition}
\label{p2.1}Let $\left( H;\left\langle \cdot ,\cdot \right\rangle \right) $
be an inner product space over the real or complex number field $\mathbb{K}$%
. The subsequent statements are equivalent.

\begin{enumerate}
\item[(i)] The following inequality holds%
\begin{equation}
\left\Vert \frac{x}{\left\Vert x\right\Vert }-\frac{y}{\left\Vert
y\right\Vert }\right\Vert \leq \left( \geq \right) r;  \label{2.1}
\end{equation}

\item[(ii)] The following reverse (improvement) of Schwarz's inequality holds%
\begin{equation}
\left\Vert x\right\Vert \left\Vert y\right\Vert -\func{Re}\left\langle
x,y\right\rangle \leq \left( \geq \right) \frac{1}{2}r^{2}\left\Vert
x\right\Vert \left\Vert y\right\Vert .  \label{2.2}
\end{equation}%
The constant $\frac{1}{2}$ is best possible in (\ref{2.2}) in the sense that
it cannot be replaced by a larger (smaller) quantity.
\end{enumerate}
\end{proposition}

\begin{proof}
It is obvious by taking the square in (\ref{2.2}) and performing the
required calculations.
\end{proof}

\begin{remark}
\label{r2.2}Since%
\begin{align*}
\left\Vert \left\Vert y\right\Vert x-\left\Vert x\right\Vert y\right\Vert &
=\left\Vert \left\Vert y\right\Vert \left( x-y\right) +\left( \left\Vert
y\right\Vert -\left\Vert x\right\Vert \right) y\right\Vert \\
& \leq \left\Vert y\right\Vert \left\Vert x-y\right\Vert +\left\vert
\left\Vert y\right\Vert -\left\Vert x\right\Vert \right\vert \left\Vert
y\right\Vert \\
& \leq 2\left\Vert y\right\Vert \left\Vert x-y\right\Vert
\end{align*}%
hence a sufficient condition for (\ref{2.1}) to hold is%
\begin{equation}
\left\Vert x-y\right\Vert \leq \frac{r}{2}\left\Vert x\right\Vert .
\label{2.3}
\end{equation}
\end{remark}

\begin{remark}
\label{r2.3}Utilising the Dunkl-Williams inequality \cite{DW}%
\begin{equation}
\left\Vert a-b\right\Vert \geq \frac{1}{2}\left( \left\Vert a\right\Vert
+\left\Vert b\right\Vert \right) \left\Vert \frac{a}{\left\Vert a\right\Vert 
}-\frac{b}{\left\Vert b\right\Vert }\right\Vert ,\ \ a,b\in H\backslash
\left\{ 0\right\}  \label{2.4}
\end{equation}%
with equality if and only if either $\left\Vert a\right\Vert =\left\Vert
b\right\Vert $ or $\left\Vert a\right\Vert +\left\Vert b\right\Vert
=\left\Vert a-b\right\Vert ,$ we can state the following inequality%
\begin{equation}
\frac{\left\Vert x\right\Vert \left\Vert y\right\Vert -\func{Re}\left\langle
x,y\right\rangle }{\left\Vert x\right\Vert \left\Vert y\right\Vert }\leq
2\left( \frac{\left\Vert x-y\right\Vert }{\left\Vert x\right\Vert
+\left\Vert y\right\Vert }\right) ^{2},\ \ x,y\in H\backslash \left\{
0\right\} .  \label{2.5}
\end{equation}%
Obviously, if $x,y\in H\backslash \left\{ 0\right\} $ are such that%
\begin{equation}
\left\Vert x-y\right\Vert \leq \eta \left( \left\Vert x\right\Vert
+\left\Vert y\right\Vert \right) ,  \label{2.6}
\end{equation}%
with $\eta \in (0,1],$ then one has the following reverse of the Schwarz
inequality%
\begin{equation}
\left\Vert x\right\Vert \left\Vert y\right\Vert -\func{Re}\left\langle
x,y\right\rangle \leq 2\eta ^{2}\left\Vert x\right\Vert \left\Vert
y\right\Vert  \label{2.7}
\end{equation}%
that is similar to (\ref{2.2}).
\end{remark}

The following result may be stated as well.

\begin{proposition}
\label{t2.4}If $x,y\in H\backslash \left\{ 0\right\} $ and $\rho >0$ are
such that%
\begin{equation}
\left\Vert \frac{x}{\left\Vert y\right\Vert }-\frac{y}{\left\Vert
x\right\Vert }\right\Vert \leq \rho ,  \label{2.8}
\end{equation}%
then we have the following reverse of Schwarz's inequality%
\begin{align}
\left( 0\leq \right) \left\Vert x\right\Vert \left\Vert y\right\Vert
-\left\vert \left\langle x,y\right\rangle \right\vert & \leq \left\Vert
x\right\Vert \left\Vert y\right\Vert -\func{Re}\left\langle x,y\right\rangle
\label{2.9} \\
& \leq \frac{1}{2}\rho ^{2}\left\Vert x\right\Vert \left\Vert y\right\Vert .
\notag
\end{align}%
The case of equality holds in the last inequality in (\ref{2.9}) if and only
if%
\begin{equation}
\left\Vert x\right\Vert =\left\Vert y\right\Vert \qquad \text{and}\qquad
\left\Vert x-y\right\Vert =\rho .  \label{2.10}
\end{equation}%
The constant $\frac{1}{2}$ in (\ref{2.9}) cannot be replaced by a smaller
quantity.
\end{proposition}

\begin{proof}
Taking the square in (\ref{2.8}), we get%
\begin{equation}
\frac{\left\Vert x\right\Vert ^{2}}{\left\Vert y\right\Vert ^{2}}-\frac{2%
\func{Re}\left\langle x,y\right\rangle }{\left\Vert x\right\Vert \left\Vert
y\right\Vert }+\frac{\left\Vert y\right\Vert ^{2}}{\left\Vert x\right\Vert
^{2}}\leq \rho ^{2}.  \label{2.11}
\end{equation}%
Since, obviously%
\begin{equation}
2\leq \frac{\left\Vert x\right\Vert ^{2}}{\left\Vert y\right\Vert ^{2}}+%
\frac{\left\Vert y\right\Vert ^{2}}{\left\Vert x\right\Vert ^{2}}
\label{2.12}
\end{equation}%
with equality iff $\left\Vert x\right\Vert =\left\Vert y\right\Vert ,$ hence
by (\ref{2.11}) we deduce the second inequality in (\ref{2.9}).

The case of equality and the best constant are obvious and we omit the
details.
\end{proof}

\begin{remark}
\label{r2.5}In \cite{H}, Hile obtained the following inequality%
\begin{equation}
\left\Vert \left\Vert x\right\Vert ^{v}x-\left\Vert y\right\Vert
^{v}y\right\Vert \leq \frac{\left\Vert x\right\Vert ^{v+1}-\left\Vert
y\right\Vert ^{v+1}}{\left\Vert x\right\Vert -\left\Vert y\right\Vert }%
\left\Vert x-y\right\Vert  \label{2.12a}
\end{equation}%
provided $v>0$ and $\left\Vert x\right\Vert \neq \left\Vert y\right\Vert .$

If in (\ref{2.12a}) we choose $v=1$ and take the square, then we get%
\begin{equation}
\left\Vert x\right\Vert ^{4}-2\left\Vert x\right\Vert \left\Vert
y\right\Vert \func{Re}\left\langle x,y\right\rangle +\left\Vert y\right\Vert
^{4}\leq \left( \left\Vert x\right\Vert +\left\Vert y\right\Vert \right)
^{2}\left\Vert x-y\right\Vert ^{2}.  \label{2.13}
\end{equation}%
Since,%
\begin{equation*}
\left\Vert x\right\Vert ^{4}+\left\Vert y\right\Vert ^{4}\geq 2\left\Vert
x\right\Vert ^{2}\left\Vert y\right\Vert ^{2},
\end{equation*}%
hence, by (\ref{2.13}) we deduce%
\begin{equation}
\left( 0\leq \right) \left\Vert x\right\Vert \left\Vert y\right\Vert -\func{%
Re}\left\langle x,y\right\rangle \leq \frac{1}{2}\cdot \frac{\left(
\left\Vert x\right\Vert +\left\Vert y\right\Vert \right) ^{2}\left\Vert
x-y\right\Vert ^{2}}{\left\Vert x\right\Vert \left\Vert y\right\Vert },
\label{2.13a}
\end{equation}%
provided $x,y\in H\backslash \left\{ 0\right\} .$
\end{remark}

The following inequality is due to Goldstein, Ryff and Clarke \cite[p. 309]%
{GRC}:%
\begin{multline}
\left\Vert x\right\Vert ^{2r}+\left\Vert y\right\Vert ^{2r}-2\left\Vert
x\right\Vert ^{r}\left\Vert y\right\Vert ^{r}\cdot \frac{\func{Re}%
\left\langle x,y\right\rangle }{\left\Vert x\right\Vert \left\Vert
y\right\Vert }  \label{2.13.a} \\
\leq \left\{ 
\begin{array}{ll}
r^{2}\left\Vert x\right\Vert ^{2r-2}\left\Vert x-y\right\Vert ^{2} & \text{%
if \ }r\geq 1 \\ 
&  \\ 
\left\Vert y\right\Vert ^{2r-2}\left\Vert x-y\right\Vert ^{2} & \text{if \ }%
r<1%
\end{array}%
\right.
\end{multline}%
provided $r\in \mathbb{R}$ and $x,y\in H$ with $\left\Vert x\right\Vert \geq
\left\Vert y\right\Vert .$

Utilising (\ref{2.13.a}) we may state the following proposition containing a
different reverse of the Schwarz inequality in inner product spaces.

\begin{proposition}
\label{p2.5.a}Let $\left( H;\left\langle \cdot ,\cdot \right\rangle \right) $
be an inner product space over the real or complex number field $\mathbb{K}$%
. If $x,y\in H\backslash \left\{ 0\right\} $ and $\left\Vert x\right\Vert
\geq \left\Vert y\right\Vert ,$ then we have%
\begin{align}
0& \leq \left\Vert x\right\Vert \left\Vert y\right\Vert -\left\vert
\left\langle x,y\right\rangle \right\vert \leq \left\Vert x\right\Vert
\left\Vert y\right\Vert -\func{Re}\left\langle x,y\right\rangle
\label{2.13.b} \\
& \leq \left\{ 
\begin{array}{ll}
\frac{1}{2}r^{2}\left( \frac{\left\Vert x\right\Vert }{\left\Vert
y\right\Vert }\right) ^{r-1}\left\Vert x-y\right\Vert ^{2} & \text{if \ }%
r\geq 1, \\ 
&  \\ 
\frac{1}{2}\left( \frac{\left\Vert x\right\Vert }{\left\Vert y\right\Vert }%
\right) ^{1-r}\left\Vert x-y\right\Vert ^{2} & \text{if \ }r<1.%
\end{array}%
\right.  \notag
\end{align}
\end{proposition}

\begin{proof}
It follows from (\ref{2.13.a}), on dividing by $\left\Vert x\right\Vert
^{r}\left\Vert y\right\Vert ^{r},$ that%
\begin{equation}
\left( \frac{\left\Vert x\right\Vert }{\left\Vert y\right\Vert }\right)
^{r}+\left( \frac{\left\Vert y\right\Vert }{\left\Vert x\right\Vert }\right)
^{r}-2\cdot \frac{\func{Re}\left\langle x,y\right\rangle }{\left\Vert
x\right\Vert \left\Vert y\right\Vert }\leq \left\{ 
\begin{array}{ll}
r^{2}\cdot \frac{\left\Vert x\right\Vert ^{r-2}}{\left\Vert y\right\Vert ^{r}%
}\left\Vert x-y\right\Vert ^{2} & \text{if \ }r\geq 1, \\ 
&  \\ 
\frac{\left\Vert y\right\Vert ^{r-2}}{\left\Vert x\right\Vert ^{r}}%
\left\Vert x-y\right\Vert ^{2} & \text{if \ }r<1.%
\end{array}%
\right.  \label{2.13.c}
\end{equation}%
Since%
\begin{equation*}
\left( \frac{\left\Vert x\right\Vert }{\left\Vert y\right\Vert }\right)
^{r}+\left( \frac{\left\Vert y\right\Vert }{\left\Vert x\right\Vert }\right)
^{r}\geq 2,
\end{equation*}%
hence, by (\ref{2.13.c}) one has%
\begin{equation*}
2-2\cdot \frac{\func{Re}\left\langle x,y\right\rangle }{\left\Vert
x\right\Vert \left\Vert y\right\Vert }\leq \left\{ 
\begin{array}{ll}
r^{2}\frac{\left\Vert x\right\Vert ^{r-2}}{\left\Vert y\right\Vert ^{r}}%
\left\Vert x-y\right\Vert ^{2} & \text{if \ }r\geq 1, \\ 
&  \\ 
\frac{\left\Vert y\right\Vert ^{r-2}}{\left\Vert x\right\Vert ^{r}}%
\left\Vert x-y\right\Vert ^{2} & \text{if \ }r<1.%
\end{array}%
\right.
\end{equation*}%
Dividing this inequality by 2 and multiplying with $\left\Vert x\right\Vert
\left\Vert y\right\Vert ,$ we deduce the desired result in (\ref{2.13.b}).
\end{proof}

Another result providing a different additive reverse (refinement) of the
Schwarz inequality may be stated.

\begin{proposition}
\label{p2.6}Let $x,y\in H$ with $y\neq 0$ and $r>0.$ The subsequent
statements are equivalent:

\begin{enumerate}
\item[(i)] The following inequality holds:%
\begin{equation}
\left\Vert x-\frac{\left\langle x,y\right\rangle }{\left\Vert y\right\Vert
^{2}}\cdot y\right\Vert \leq \left( \geq \right) r;  \label{2.14}
\end{equation}

\item[(ii)] The following reverse (refinement) of the quadratic Schwarz
inequality holds:%
\begin{equation}
\left\Vert x\right\Vert ^{2}\left\Vert y\right\Vert ^{2}-\left\vert
\left\langle x,y\right\rangle \right\vert ^{2}\leq \left( \geq \right)
r^{2}\left\Vert y\right\Vert ^{2}.  \label{2.15}
\end{equation}
\end{enumerate}
\end{proposition}

The proof is obvious on taking the square in (\ref{2.14}) and performing the
calculation.

\begin{remark}
\label{r2.7}Since%
\begin{align*}
\left\Vert \left\Vert y\right\Vert ^{2}x-\left\langle x,y\right\rangle
y\right\Vert & =\left\Vert \left\Vert y\right\Vert ^{2}\left( x-y\right)
-\left\langle x-y,y\right\rangle y\right\Vert \\
& \leq \left\Vert y\right\Vert ^{2}\left\Vert x-y\right\Vert +\left\vert
\left\langle x-y,y\right\rangle \right\vert \left\Vert y\right\Vert \\
& \leq 2\left\Vert x-y\right\Vert \left\Vert y\right\Vert ^{2},
\end{align*}%
hence a sufficient condition for the inequality (\ref{2.14}) to hold is that%
\begin{equation}
\left\Vert x-y\right\Vert \leq \frac{r}{2}.  \label{2.16}
\end{equation}
\end{remark}

The following proposition may give a complementary approach:

\begin{proposition}
\label{p2.6.a}Let $x,y\in H$ with $\left\langle x,y\right\rangle \neq 0$ and 
$\rho >0.$ If 
\begin{equation}
\left\Vert x-\frac{\left\langle x,y\right\rangle }{\left\vert \left\langle
x,y\right\rangle \right\vert }\cdot y\right\Vert \leq \rho ,  \label{2.16.a}
\end{equation}%
then%
\begin{equation}
(0\leq )\left\Vert x\right\Vert \left\Vert y\right\Vert -\left\vert
\left\langle x,y\right\rangle \right\vert \leq \frac{1}{2}\rho ^{2}.
\label{2.16.b}
\end{equation}%
The multiplicative constant $\frac{1}{2}$ is best possible in (\ref{2.16.b}).
\end{proposition}

The proof is similar to the ones outlined above and we omit it.

For the case of complex inner product spaces, we may state the following
result.

\begin{proposition}
\label{t2.8}Let $\left( H;\left\langle \cdot ,\cdot \right\rangle \right) $
be a complex inner product space and $\alpha \in \mathbb{C}$ a given complex
number with $\func{Re}\alpha ,$ $\func{Im}\alpha >0.$ If $x,y\in H$ are such
that%
\begin{equation}
\left\Vert x-\frac{\func{Im}\alpha }{\func{Re}\alpha }\cdot y\right\Vert
\leq r,  \label{2.17}
\end{equation}%
then we have the inequality%
\begin{align}
(0& \leq )\left\Vert x\right\Vert \left\Vert y\right\Vert -\left\vert
\left\langle x,y\right\rangle \right\vert \leq \left\Vert x\right\Vert
\left\Vert y\right\Vert -\func{Re}\left\langle x,y\right\rangle  \label{2.18}
\\
& \leq \frac{1}{2}\cdot \frac{\func{Re}\alpha }{\func{Im}\alpha }\cdot r^{2}.
\notag
\end{align}%
The equality holds in the second inequality in (\ref{2.18}) if and only if
the case of equality holds in (\ref{2.17}) and $\func{Re}\alpha \cdot
\left\Vert x\right\Vert =\func{Im}\alpha \cdot \left\Vert y\right\Vert .$
\end{proposition}

\begin{proof}
Observe that the condition (\ref{2.17}) is equivalent to%
\begin{equation}
\left[ \func{Re}\alpha \right] ^{2}\left\Vert x\right\Vert ^{2}+\left[ \func{%
Im}\alpha \right] ^{2}\left\Vert y\right\Vert ^{2}\leq 2\func{Re}\alpha 
\func{Im}\alpha \func{Re}\left\langle x,y\right\rangle +\left[ \func{Re}%
\alpha \right] ^{2}r^{2}.  \label{2.19}
\end{equation}%
On the other hand, on utilising the elementary inequality%
\begin{equation}
2\func{Re}\alpha \func{Im}\alpha \left\Vert x\right\Vert \left\Vert
y\right\Vert \leq \left[ \func{Re}\alpha \right] ^{2}\left\Vert x\right\Vert
^{2}+\left[ \func{Im}\alpha \right] ^{2}\left\Vert y\right\Vert ^{2},
\label{2.20}
\end{equation}%
with equality if and only if $\func{Re}\alpha \cdot \left\Vert x\right\Vert =%
\func{Im}\alpha \cdot \left\Vert y\right\Vert ,$ we deduce from (\ref{2.19})
that%
\begin{equation}
2\func{Re}\alpha \func{Im}\alpha \left\Vert x\right\Vert \left\Vert
y\right\Vert \leq 2\func{Re}\alpha \func{Im}\alpha \func{Re}\left\langle
x,y\right\rangle +r^{2}\left[ \func{Re}\alpha \right] ^{2}  \label{2.21}
\end{equation}%
giving the desired inequality (\ref{2.18}).

The case of equality follows from the above and we omit the details.
\end{proof}

The following different reverse for the Schwarz inequality that holds for
both real and complex inner product spaces may be stated as well.

\begin{theorem}
\label{t2.9}Let $\left( H;\left\langle \cdot ,\cdot \right\rangle \right) $
be an inner product space over $\mathbb{K}$, $\mathbb{K}=\mathbb{C},\mathbb{R%
}.$ If $\alpha \in \mathbb{K}\backslash \left\{ 0\right\} ,$ then%
\begin{align}
0& \leq \left\Vert x\right\Vert \left\Vert y\right\Vert -\left\vert
\left\langle x,y\right\rangle \right\vert \leq \left\Vert x\right\Vert
\left\Vert y\right\Vert -\func{Re}\left[ \frac{\alpha ^{2}}{\left\vert
\alpha \right\vert ^{2}}\left\langle x,y\right\rangle \right]  \label{2.22}
\\
& \leq \frac{1}{2}\cdot \frac{\left[ \left\vert \func{Re}\alpha \right\vert
\left\Vert x-y\right\Vert +\left\vert \func{Im}\alpha \right\vert \left\Vert
x+y\right\Vert \right] ^{2}}{\left\vert \alpha \right\vert ^{2}}\leq \frac{1%
}{2}\cdot I^{2},  \notag
\end{align}%
where%
\begin{equation}
I:=\left\{ 
\begin{array}{ll}
\max \left\{ \left\vert \func{Re}\alpha \right\vert ,\left\vert \func{Im}%
\alpha \right\vert \right\} \left( \left\Vert x-y\right\Vert +\left\Vert
x+y\right\Vert \right) ; &  \\ 
&  \\ 
\left( \left\vert \func{Re}\alpha \right\vert ^{p}+\left\vert \func{Im}%
\alpha \right\vert ^{p}\right) ^{\frac{1}{p}}\left( \left\Vert
x-y\right\Vert ^{q}+\left\Vert x+y\right\Vert ^{q}\right) ^{\frac{1}{q}}, & 
p>1, \\ 
& \frac{1}{p}+\frac{1}{q}=1; \\ 
\max \left\{ \left\Vert x-y\right\Vert ,\left\Vert x+y\right\Vert \right\}
\left( \left\vert \func{Re}\alpha \right\vert +\left\vert \func{Im}\alpha
\right\vert \right) . & 
\end{array}%
\right.  \label{2.23}
\end{equation}
\end{theorem}

\begin{proof}
Observe, for $\alpha \in \mathbb{K}\backslash \left\{ 0\right\} ,$ that%
\begin{align*}
\left\Vert \alpha x-\bar{\alpha}y\right\Vert ^{2}& =\left\vert \alpha
\right\vert ^{2}\left\Vert x\right\Vert ^{2}-2\func{Re}\left\langle \alpha x,%
\bar{\alpha}y\right\rangle +\left\vert \alpha \right\vert ^{2}\left\Vert
y\right\Vert ^{2} \\
& =\left\vert \alpha \right\vert ^{2}\left( \left\Vert x\right\Vert
^{2}+\left\Vert y\right\Vert ^{2}\right) -2\func{Re}\left[ \alpha
^{2}\left\langle x,y\right\rangle \right] .
\end{align*}%
Since $\left\Vert x\right\Vert ^{2}+\left\Vert y\right\Vert ^{2}\geq
2\left\Vert x\right\Vert \left\Vert y\right\Vert ,$ hence%
\begin{equation}
\left\Vert \alpha x-\bar{\alpha}y\right\Vert ^{2}\geq 2\left\vert \alpha
\right\vert ^{2}\left\{ \left\Vert x\right\Vert \left\Vert y\right\Vert -%
\func{Re}\left[ \frac{\alpha ^{2}}{\left\vert \alpha \right\vert ^{2}}%
\left\langle x,y\right\rangle \right] \right\} .  \label{2.24}
\end{equation}%
On the other hand, we have%
\begin{align}
\left\Vert \alpha x-\bar{\alpha}y\right\Vert & =\left\Vert \left( \func{Re}%
\alpha +i\func{Im}\alpha \right) x-\left( \func{Re}\alpha -i\func{Im}\alpha
\right) y\right\Vert   \label{2.25} \\
& =\left\Vert \func{Re}\alpha \left( x-y\right) +i\func{Im}\alpha \left(
x+y\right) \right\Vert   \notag \\
& \leq \left\vert \func{Re}\alpha \right\vert \left\Vert x-y\right\Vert
+\left\vert \func{Im}\alpha \right\vert \left\Vert x+y\right\Vert .  \notag
\end{align}%
Utilising (\ref{2.24}) and (\ref{2.25}) we deduce the third inequality in (%
\ref{2.22}).

For the last inequality we use the following elementary inequality%
\begin{equation}
\alpha a+\beta b\leq \left\{ 
\begin{array}{ll}
\max \left\{ \alpha ,\beta \right\} \left( a+b\right) &  \\ 
&  \\ 
\left( \alpha ^{p}+\beta ^{p}\right) ^{\frac{1}{p}}\left( a^{q}+b^{q}\right)
^{\frac{1}{q}}, & p>1,\ \frac{1}{p}+\frac{1}{q}=1,%
\end{array}%
\right.  \label{2.26}
\end{equation}%
provided $\alpha ,\beta ,a,b\geq 0.$
\end{proof}

The following result may be stated.

\begin{proposition}
\label{p2.11}Let $\left( H;\left\langle \cdot ,\cdot \right\rangle \right) $
be an inner product over$\mathbb{\ K}$ and $e\in H,$ $\left\Vert
e\right\Vert =1.$ If $\lambda \in \left( 0,1\right) ,$ then%
\begin{multline}
\func{Re}\left[ \left\langle x,y\right\rangle -\left\langle x,e\right\rangle
\left\langle e,y\right\rangle \right]  \label{2.28} \\
\leq \frac{1}{4}\cdot \frac{1}{\lambda \left( 1-\lambda \right) }\left[
\left\Vert \lambda x+\left( 1-\lambda \right) y\right\Vert ^{2}-\left\vert
\left\langle \lambda x+\left( 1-\lambda \right) y,e\right\rangle \right\vert
^{2}\right] .
\end{multline}%
The constant $\frac{1}{4}$ is best possible.
\end{proposition}

\begin{proof}
Firstly, note that the following equality holds true%
\begin{equation*}
\left\langle x-\left\langle x,e\right\rangle e,y-\left\langle
y,e\right\rangle e\right\rangle =\left\langle x,y\right\rangle -\left\langle
x,e\right\rangle \left\langle e,y\right\rangle .
\end{equation*}%
Utilising the elementary inequality%
\begin{equation*}
\func{Re}\left\langle z,w\right\rangle \leq \frac{1}{4}\left\Vert
z+w\right\Vert ^{2},\qquad z,w\in H
\end{equation*}%
we have%
\begin{align*}
& \func{Re}\left\langle x-\left\langle x,e\right\rangle e,y-\left\langle
y,e\right\rangle e\right\rangle \\
& =\frac{1}{\lambda \left( 1-\lambda \right) }\func{Re}\left\langle \lambda
x-\left\langle \lambda x,e\right\rangle e,\left( 1-\lambda \right)
y-\left\langle \left( 1-\lambda \right) y,e\right\rangle e\right\rangle \\
& \leq \frac{1}{4}\cdot \frac{1}{\lambda \left( 1-\lambda \right) }\left[
\left\Vert \lambda x+\left( 1-\lambda \right) y\right\Vert ^{2}-\left\vert
\left\langle \lambda x+\left( 1-\lambda \right) y,e\right\rangle \right\vert
^{2}\right] ,
\end{align*}%
proving the desired inequality (\ref{2.28}).
\end{proof}

\begin{remark}
\label{r2.12}For $\lambda =\frac{1}{2},$ we get the simpler inequality:%
\begin{equation}
\func{Re}\left[ \left\langle x,y\right\rangle -\left\langle x,e\right\rangle
\left\langle e,y\right\rangle \right] \leq \left\Vert \frac{x+y}{2}%
\right\Vert ^{2}-\left\vert \left\langle \frac{x+y}{2},e\right\rangle
\right\vert ^{2},  \label{2.29}
\end{equation}%
that has been obtained in \cite[p. 46]{SSD3}, for which the sharpness of the
inequality was established.
\end{remark}

The following result may be stated as well.

\begin{theorem}
\label{t2.13}Let $\left( H;\left\langle \cdot ,\cdot \right\rangle \right) $
be an inner product space over $\mathbb{K}$ and $p\geq 1.$ Then for any $%
x,y\in H$ we have%
\begin{align}
0& \leq \left\Vert x\right\Vert \left\Vert y\right\Vert -\left\vert
\left\langle x,y\right\rangle \right\vert \leq \left\Vert x\right\Vert
\left\Vert y\right\Vert -\func{Re}\left\langle x,y\right\rangle  \label{2.30}
\\
& \leq \frac{1}{2}\times \left\{ 
\begin{array}{l}
\left[ \left( \left\Vert x\right\Vert +\left\Vert y\right\Vert \right)
^{2p}-\left\Vert x+y\right\Vert ^{2p}\right] ^{\frac{1}{p}}, \\ 
\\ 
\left[ \left\Vert x-y\right\Vert ^{2p}-\left\vert \left\Vert x\right\Vert
-\left\Vert y\right\Vert \right\vert ^{2p}\right] ^{\frac{1}{p}}.%
\end{array}%
\right.  \notag
\end{align}
\end{theorem}

\begin{proof}
Firstly, observe that%
\begin{equation*}
2\left( \left\Vert x\right\Vert \left\Vert y\right\Vert -\func{Re}%
\left\langle x,y\right\rangle \right) =\left( \left\Vert x\right\Vert
+\left\Vert y\right\Vert \right) ^{2}-\left\Vert x+y\right\Vert ^{2}.
\end{equation*}%
Denoting $D:=\left\Vert x\right\Vert \left\Vert y\right\Vert -\func{Re}%
\left\langle x,y\right\rangle ,$ then we have%
\begin{equation}
2D+\left\Vert x+y\right\Vert ^{2}=\left( \left\Vert x\right\Vert +\left\Vert
y\right\Vert \right) ^{2}.  \label{2.31}
\end{equation}%
Taking in (\ref{2.31}) the power $p\geq 1$ and using the elementary
inequality 
\begin{equation}
\left( a+b\right) ^{p}\geq a^{p}+b^{p};a,b\geq 0,  \label{2.31.a}
\end{equation}%
we have%
\begin{equation*}
\left( \left\Vert x\right\Vert +\left\Vert y\right\Vert \right) ^{2p}=\left(
2D+\left\Vert x+y\right\Vert ^{2}\right) ^{p}\geq 2^{p}D^{p}+\left\Vert
x+y\right\Vert ^{2p}
\end{equation*}%
giving%
\begin{equation*}
D^{p}\leq \frac{1}{2^{p}}\left[ \left( \left\Vert x\right\Vert +\left\Vert
y\right\Vert \right) ^{2p}-\left\Vert x+y\right\Vert ^{2p}\right] ,
\end{equation*}%
which is clearly equivalent to the first branch of the third inequality in (%
\ref{2.30}).

With the above notation, we also have%
\begin{equation}
2D+\left( \left\Vert x\right\Vert -\left\Vert y\right\Vert \right)
^{2}=\left\Vert x-y\right\Vert ^{2}.  \label{2.32}
\end{equation}%
Taking the power $p\geq 1$ in (\ref{2.32}) and using the inequality (\ref%
{2.31.a}) we deduce%
\begin{equation*}
\left\Vert x-y\right\Vert ^{2p}\geq 2^{p}D^{p}+\left\vert \left\Vert
x\right\Vert -\left\Vert y\right\Vert \right\vert ^{2p},
\end{equation*}%
from where we get the last part of (\ref{2.30}).
\end{proof}

\section{More Schwarz's Related Inequalities}

Before we point out other inequalities related to the Schwarz inequality, we
need the following identity that is interesting in itself.

\begin{lemma}
\label{l2.14}Let $\left( H;\left\langle \cdot ,\cdot \right\rangle \right) $
be an inner product space over the real or complex number field $\mathbb{K}$%
, $e\in H,$ $\left\Vert e\right\Vert =1,$ $\alpha \in H$ and $\gamma ,\Gamma
\in \mathbb{K}$. Then we have the identity:%
\begin{multline}
\left\Vert x\right\Vert ^{2}-\left\vert \left\langle x,e\right\rangle
\right\vert ^{2}+\left\Vert x-\frac{\gamma +\Gamma }{2}e\right\Vert ^{2}
\label{2.33} \\
=\left( \func{Re}\Gamma -\func{Re}\left\langle x,e\right\rangle \right)
\left( \func{Re}\left\langle x,e\right\rangle -\func{Re}\gamma \right) \\
+\left( \func{Im}\Gamma -\func{Im}\left\langle x,e\right\rangle \right)
\left( \func{Im}\left\langle x,e\right\rangle -\func{Im}\gamma \right) +%
\frac{1}{4}\left\vert \Gamma -\gamma \right\vert ^{2}.
\end{multline}
\end{lemma}

\begin{proof}
We start with the following known equality (see for instance \cite[eq. (2.6)]%
{SSD1a})%
\begin{equation}
\left\Vert x\right\Vert ^{2}-\left\vert \left\langle x,e\right\rangle
\right\vert ^{2}=\func{Re}\left[ \left( \Gamma -\left\langle
x,e\right\rangle \right) \left( \overline{\left\langle x,e\right\rangle }-%
\bar{\gamma}\right) \right] -\func{Re}\left\langle \Gamma e-x,x-\gamma
e\right\rangle  \label{2.34}
\end{equation}%
holding for $x\in H,$ $e\in H,$ $\left\Vert e\right\Vert =1$ and $\gamma
,\Gamma \in \mathbb{K}$.

We also know that (see for instance \cite{SSD1b})%
\begin{equation}
\func{Re}\left\langle \Gamma e-x,x-\gamma e\right\rangle =\left\Vert x-\frac{%
\gamma +\Gamma }{2}e\right\Vert ^{2}-\frac{1}{4}\left\vert \Gamma -\gamma
\right\vert ^{2}.  \label{2.35}
\end{equation}%
Since%
\begin{align}
& \func{Re}\left[ \left( \Gamma -\left\langle x,e\right\rangle \right)
\left( \overline{\left\langle x,e\right\rangle }-\bar{\gamma}\right) \right]
\label{2.36} \\
& =\left( \func{Re}\Gamma -\func{Re}\left\langle x,e\right\rangle \right)
\left( \func{Re}\left\langle x,e\right\rangle -\func{Re}\gamma \right) 
\notag \\
& \qquad \qquad +\left( \func{Im}\Gamma -\func{Im}\left\langle
x,e\right\rangle \right) \left( \func{Im}\left\langle x,e\right\rangle -%
\func{Im}\gamma \right) ,  \notag
\end{align}%
hence, by (\ref{2.34}) -- (\ref{2.36}), we deduce the desired identity (\ref%
{2.33}).
\end{proof}

The following general result providing a reverse of the Schwarz inequality
may be stated.

\begin{proposition}
\label{t2.15}Let $\left( H;\left\langle \cdot ,\cdot \right\rangle \right) $
be an inner product space over $\mathbb{K},$ $e\in H,$ $\left\Vert
e\right\Vert =1,$ $x\in H$ and $\gamma ,\Gamma \in \mathbb{K}$. Then we have
the inequality:%
\begin{equation}
\left( 0\leq \right) \left\Vert x\right\Vert ^{2}-\left\vert \left\langle
x,e\right\rangle \right\vert ^{2}+\left\Vert x-\frac{\gamma +\Gamma }{2}%
\cdot e\right\Vert ^{2}\leq \frac{1}{2}\left\vert \Gamma -\gamma \right\vert
^{2}.  \label{2.37}
\end{equation}%
The constant $\frac{1}{2}$ is best possible in (\ref{2.37}). The case of
equality holds in (\ref{2.37}) if and only if%
\begin{equation}
\func{Re}\left\langle x,e\right\rangle =\func{Re}\left( \frac{\gamma +\Gamma 
}{2}\right) ,\qquad \func{Im}\left\langle x,e\right\rangle =\func{Im}\left( 
\frac{\gamma +\Gamma }{2}\right) .  \label{2.38}
\end{equation}
\end{proposition}

\begin{proof}
Utilising the elementary inequality for real numbers%
\begin{equation*}
\alpha \beta \leq \frac{1}{4}\left( \alpha +\beta \right) ^{2},\qquad \alpha
,\beta \in \mathbb{R};
\end{equation*}%
with equality iff $\alpha =\beta ,$ we have%
\begin{equation}
\left( \func{Re}\Gamma -\func{Re}\left\langle x,e\right\rangle \right)
\left( \func{Re}\left\langle x,e\right\rangle -\func{Re}\gamma \right) \leq 
\frac{1}{4}\left( \func{Re}\Gamma -\func{Re}\gamma \right) ^{2}  \label{2.39}
\end{equation}%
and%
\begin{equation}
\left( \func{Im}\Gamma -\func{Im}\left\langle x,e\right\rangle \right)
\left( \func{Im}\left\langle x,e\right\rangle -\func{Im}\gamma \right) \leq 
\frac{1}{4}\left( \func{Im}\Gamma -\func{Im}\gamma \right) ^{2}  \label{2.40}
\end{equation}%
with equality if and only if%
\begin{equation*}
\func{Re}\left\langle x,e\right\rangle =\frac{\func{Re}\Gamma +\func{Re}%
\gamma }{2}\qquad \text{and\qquad }\func{Im}\left\langle x,e\right\rangle =%
\frac{\func{Im}\Gamma +\func{Im}\gamma }{2}.
\end{equation*}%
Finally, on making use of (\ref{2.39}), (\ref{2.40}) and the identity (\ref%
{2.33}), we deduce the desired result (\ref{2.37}).
\end{proof}

The following result may be stated as well.

\begin{proposition}
\label{t2.16}Let $\left( H;\left\langle \cdot ,\cdot \right\rangle \right) $
be an inner product space over $\mathbb{K},$ $e\in H,$ $\left\Vert
e\right\Vert =1,$ $x\in H$ and $\gamma ,\Gamma \in \mathbb{K}$. If $x\in H$
is such that%
\begin{equation}
\func{Re}\gamma \leq \func{Re}\left\langle x,e\right\rangle \leq \func{Re}%
\Gamma \qquad \text{and\qquad }\func{Im}\gamma \leq \func{Im}\left\langle
x,e\right\rangle \leq \func{Im}\Gamma ,  \label{2.41}
\end{equation}%
then we have the inequality%
\begin{equation}
\left\Vert x\right\Vert ^{2}-\left\vert \left\langle x,e\right\rangle
\right\vert ^{2}+\left\Vert x-\frac{\gamma +\Gamma }{2}e\right\Vert ^{2}\geq 
\frac{1}{4}\left\vert \Gamma -\gamma \right\vert ^{2}.  \label{2.42}
\end{equation}%
The constant $\frac{1}{4}$ is best possible in (\ref{2.42}). The case of
equality holds in (\ref{2.42}) if and only if%
\begin{equation*}
\func{Re}\left\langle x,e\right\rangle =\func{Re}\Gamma \text{ or }\func{Re}%
\left\langle x,e\right\rangle =\func{Re}\gamma 
\end{equation*}%
and%
\begin{equation*}
\func{Im}\left\langle x,e\right\rangle =\func{Im}\Gamma \text{ or }\func{Im}%
\left\langle x,e\right\rangle =\func{Im}\gamma .
\end{equation*}
\end{proposition}

\begin{proof}
From the hypothesis we obviously have%
\begin{equation*}
\left( \func{Re}\Gamma -\func{Re}\left\langle x,e\right\rangle \right)
\left( \func{Re}\left\langle x,e\right\rangle -\func{Re}\gamma \right) \geq 0
\end{equation*}%
and%
\begin{equation*}
\left( \func{Im}\Gamma -\func{Im}\left\langle x,e\right\rangle \right)
\left( \func{Im}\left\langle x,e\right\rangle -\func{Im}\gamma \right) \geq
0.
\end{equation*}%
Utilising the identity (\ref{2.33}) we deduce the desired result (\ref{2.42}%
). The case of equality is obvious.
\end{proof}

Further on, we can state the following reverse of the quadratic Schwarz
inequality:

\begin{proposition}
\label{t2.17}Let $\left( H;\left\langle \cdot ,\cdot \right\rangle \right) $
be an inner product space over $\mathbb{K},$ $e\in H,$ $\left\Vert
e\right\Vert =1.$ If $\gamma ,\Gamma \in \mathbb{K}$ and $x\in H$ are such
that either%
\begin{equation}
\func{Re}\left\langle \Gamma e-x,x-\gamma e\right\rangle \geq 0  \label{2.43}
\end{equation}%
or, equivalently,%
\begin{equation}
\left\Vert x-\frac{\gamma +\Gamma }{2}e\right\Vert \leq \frac{1}{2}%
\left\vert \Gamma -\gamma \right\vert ,  \label{2.44}
\end{equation}%
then%
\begin{align}
(0& \leq )\left\Vert x\right\Vert ^{2}-\left\vert \left\langle
x,e\right\rangle \right\vert ^{2}  \label{2.45} \\
& \leq \left( \func{Re}\Gamma -\func{Re}\left\langle x,e\right\rangle
\right) \left( \func{Re}\left\langle x,e\right\rangle -\func{Re}\gamma
\right)  \notag \\
& \qquad \qquad \qquad +\left( \func{Im}\Gamma -\func{Im}\left\langle
x,e\right\rangle \right) \left( \func{Im}\left\langle x,e\right\rangle -%
\func{Im}\gamma \right)  \notag \\
& \leq \frac{1}{4}\left\vert \Gamma -\gamma \right\vert ^{2}.  \notag
\end{align}%
The case of equality holds in (\ref{2.45}) if it holds either in (\ref{2.43}%
) or (\ref{2.44}).
\end{proposition}

The proof is obvious by Lemma \ref{l2.14} and we omit the details.

\begin{remark}
\label{r2.18}We remark that the inequality (\ref{2.45}) may also be used to
get, for instance, the following result%
\begin{multline}
\left\Vert x\right\Vert ^{2}-\left\vert \left\langle x,e\right\rangle
\right\vert ^{2}\leq \left[ \left( \func{Re}\Gamma -\func{Re}\left\langle
x,e\right\rangle \right) ^{2}+\left( \func{Im}\Gamma -\func{Im}\left\langle
x,e\right\rangle \right) ^{2}\right] ^{\frac{1}{2}}  \label{2.46} \\
\times \left[ \left( \func{Re}\left\langle x,e\right\rangle -\func{Re}\gamma
\right) ^{2}+\left( \func{Im}\left\langle x,e\right\rangle -\func{Im}\gamma
\right) ^{2}\right] ^{\frac{1}{2}},
\end{multline}%
that provides a different bound than $\frac{1}{4}\left\vert \Gamma -\gamma
\right\vert ^{2}$ for the quantity $\left\Vert x\right\Vert ^{2}-\left\vert
\left\langle x,e\right\rangle \right\vert ^{2}.$
\end{remark}

The following result may be stated as well.

\begin{theorem}
\label{t2.19}Let $\left( H;\left\langle \cdot ,\cdot \right\rangle \right) $
be an inner product space over $\mathbb{K}$ and $\alpha ,\gamma >0,$ $\beta
\in \mathbb{K}$ with $\left\vert \beta \right\vert ^{2}\geq \alpha \gamma .$
If $x,a\in H$ are such that $a\neq 0$ and%
\begin{equation}
\left\Vert x-\frac{\beta }{\alpha }a\right\Vert \leq \frac{\left( \left\vert
\beta \right\vert ^{2}-\alpha \gamma \right) ^{\frac{1}{2}}}{\alpha }%
\left\Vert a\right\Vert ,  \label{2.47}
\end{equation}%
then we have the following reverses of Schwarz's inequality%
\begin{align}
\left\Vert x\right\Vert \left\Vert a\right\Vert & \leq \frac{\func{Re}\beta
\cdot \func{Re}\left\langle x,a\right\rangle +\func{Im}\beta \cdot \func{Im}%
\left\langle x,a\right\rangle }{\sqrt{\alpha \gamma }}  \label{2.48} \\
& \leq \frac{\left\vert \beta \right\vert \left\vert \left\langle
x,a\right\rangle \right\vert }{\sqrt{\alpha \gamma }}  \notag
\end{align}%
and%
\begin{equation}
\left( 0\leq \right) \left\Vert x\right\Vert ^{2}\left\Vert a\right\Vert
^{2}-\left\vert \left\langle x,a\right\rangle \right\vert ^{2}\leq \frac{%
\left\vert \beta \right\vert ^{2}-\alpha \gamma }{\alpha \gamma }\left\vert
\left\langle x,a\right\rangle \right\vert ^{2}.  \label{2.49}
\end{equation}
\end{theorem}

\begin{proof}
Taking the square in (\ref{2.47}), it becomes equivalent to%
\begin{equation*}
\left\Vert x\right\Vert ^{2}-\frac{2}{\alpha }\func{Re}\left[ \bar{\beta}%
\left\langle x,a\right\rangle \right] +\frac{\left\vert \beta \right\vert
^{2}}{\alpha ^{2}}\left\Vert a\right\Vert ^{2}\leq \frac{\left\vert \beta
\right\vert ^{2}-\alpha \gamma }{\alpha ^{2}}\left\Vert a\right\Vert ^{2},
\end{equation*}%
which is clearly equivalent to%
\begin{align}
\alpha \left\Vert x\right\Vert ^{2}+\gamma \left\Vert a\right\Vert ^{2}&
\leq 2\func{Re}\left[ \bar{\beta}\left\langle x,a\right\rangle \right]
\label{2.50} \\
& =2\left[ \func{Re}\beta \cdot \func{Re}\left\langle x,a\right\rangle +%
\func{Im}\beta \cdot \func{Im}\left\langle x,a\right\rangle \right] .  \notag
\end{align}%
On the other hand, since%
\begin{equation}
2\sqrt{\alpha \gamma }\left\Vert x\right\Vert \left\Vert a\right\Vert \leq
\alpha \left\Vert x\right\Vert ^{2}+\gamma \left\Vert a\right\Vert ^{2},
\label{2.51}
\end{equation}%
hence by (\ref{2.50}) and (\ref{2.51}) we deduce the first inequality in (%
\ref{2.48}).

The other inequalities are obvious.
\end{proof}

\begin{remark}
\label{r2.20}The above inequality (\ref{2.48}) contains in particular the
reverse (\ref{1.11}) of the Schwarz inequality. Indeed, if we assume that $%
\alpha =1,$ $\beta =\frac{\delta +\Delta }{2},$ $\delta ,\Delta \in \mathbb{K%
}$, with $\gamma =\func{Re}\left( \Delta \bar{\gamma}\right) >0,$ then the
condition $\left\vert \beta \right\vert ^{2}\geq \alpha \gamma $ is
equivalent to $\left\vert \delta +\Delta \right\vert ^{2}\geq 4\func{Re}%
\left( \Delta \bar{\gamma}\right) $ which is actually $\left\vert \Delta
-\delta \right\vert ^{2}\geq 0.$ With this assumption, (\ref{2.47}) becomes%
\begin{equation*}
\left\Vert x-\frac{\delta +\Delta }{2}\cdot a\right\Vert \leq \frac{1}{2}%
\left\vert \Delta -\delta \right\vert \left\Vert a\right\Vert ,
\end{equation*}%
which implies the reverse of the Schwarz inequality%
\begin{align*}
\left\Vert x\right\Vert \left\Vert a\right\Vert & \leq \frac{\func{Re}\left[
\left( \bar{\Delta}+\bar{\delta}\right) \left\langle x,a\right\rangle \right]
}{2\sqrt{\func{Re}\left( \Delta \bar{\delta}\right) }} \\
& \leq \frac{\left\vert \Delta +\delta \right\vert }{2\sqrt{\func{Re}\left(
\Delta \bar{\delta}\right) }}\left\vert \left\langle x,a\right\rangle
\right\vert ,
\end{align*}%
which is (\ref{1.11}).
\end{remark}

The following particular case of Theorem \ref{t2.19} may be stated:

\begin{corollary}
\label{c2.21}Let $\left( H;\left\langle \cdot ,\cdot \right\rangle \right) $
be an inner product space over $\mathbb{K},$ $\varphi \in \lbrack 0,2\pi ),$ 
$\theta \in \left( 0,\frac{\pi }{2}\right) .$ If $x,a\in H$ are such that $%
a\neq 0$ and%
\begin{equation}
\left\Vert x-\left( \cos \varphi +i\sin \varphi \right) a\right\Vert \leq
\cos \theta \left\Vert a\right\Vert ,  \label{2.52}
\end{equation}%
then we have the reverses of the Schwarz inequality%
\begin{equation}
\left\Vert x\right\Vert \left\Vert a\right\Vert \leq \frac{\cos \varphi 
\func{Re}\left\langle x,a\right\rangle +\sin \varphi \func{Im}\left\langle
x,a\right\rangle }{\sin \theta }.  \label{2.52a}
\end{equation}%
In particular, if%
\begin{equation*}
\left\Vert x-a\right\Vert \leq \cos \theta \left\Vert a\right\Vert ,
\end{equation*}%
then%
\begin{equation*}
\left\Vert x\right\Vert \left\Vert a\right\Vert \leq \frac{1}{\cos \theta }%
\func{Re}\left\langle x,a\right\rangle ;
\end{equation*}%
and if%
\begin{equation*}
\left\Vert x-ia\right\Vert \leq \cos \theta \left\Vert a\right\Vert ,
\end{equation*}%
then%
\begin{equation*}
\left\Vert x\right\Vert \left\Vert a\right\Vert \leq \frac{1}{\cos \theta }%
\func{Im}\left\langle x,a\right\rangle .
\end{equation*}
\end{corollary}

\section{Reverses of the Generalised Triangle Inequality}

In \cite{SSD4}, the author obtained the following reverse result for the
generalised triangle inequality%
\begin{equation}
\sum_{i=1}^{n}\left\Vert x_{i}\right\Vert \geq \left\Vert
\sum_{i=1}^{n}x_{i}\right\Vert ,  \label{3.1}
\end{equation}%
provided $x_{i}\in H,$ $i\in \left\{ 1,\dots ,n\right\} $ are vectors in a
real or complex inner product $\left( H;\left\langle \cdot ,\cdot
\right\rangle \right) :$

\begin{theorem}
\label{t3.1}Let $e,x_{i}\in H,$ $i\in \left\{ 1,\dots ,n\right\} $ with $%
\left\Vert e\right\Vert =1.$ If $k_{i}\geq 0,$ $i\in \left\{ 1,\dots
,n\right\} $ are such that%
\begin{equation}
\left( 0\leq \right) \left\Vert x_{i}\right\Vert -\func{Re}\left\langle
e,x_{i}\right\rangle \leq k_{i}\qquad \text{for each \qquad }i\in \left\{
1,\dots ,n\right\} ,  \label{3.2}
\end{equation}%
then we have the inequality%
\begin{equation}
\left( 0\leq \right) \sum_{i=1}^{n}\left\Vert x_{i}\right\Vert -\left\Vert
\sum_{i=1}^{n}x_{i}\right\Vert \leq \sum_{i=1}^{n}k_{i}.  \label{3.3}
\end{equation}%
The equality holds in (\ref{3.3}) if and only if%
\begin{equation}
\sum_{i=1}^{n}\left\Vert x_{i}\right\Vert \geq \sum_{i=1}^{n}k_{i}
\label{3.4}
\end{equation}%
and%
\begin{equation}
\sum_{i=1}^{n}x_{i}=\left( \sum_{i=1}^{n}\left\Vert x_{i}\right\Vert
-\sum_{i=1}^{n}k_{i}\right) e.  \label{3.5}
\end{equation}
\end{theorem}

By utilising some of the results obtained in Section \ref{s2}, we point out
several reverses of the generalised triangle inequality (\ref{3.1}) that are
corollaries of the above Theorem \ref{t3.1}.

\begin{corollary}
\label{c3.2}Let $e,$ $x_{i}\in H\backslash \left\{ 0\right\} ,$ $i\in
\left\{ 1,\dots ,n\right\} $ with $\left\Vert e\right\Vert =1.$ If%
\begin{equation}
\left\Vert \frac{x_{i}}{\left\Vert x_{i}\right\Vert }-e\right\Vert \leq
r_{i}\qquad \text{for each \qquad }i\in \left\{ 1,\dots ,n\right\} ,
\label{3.6}
\end{equation}%
then%
\begin{align}
(0& \leq )\sum_{i=1}^{n}\left\Vert x_{i}\right\Vert -\left\Vert
\sum_{i=1}^{n}x_{i}\right\Vert  \label{3.7} \\
& \leq \frac{1}{2}\sum_{i=1}^{n}r_{i}^{2}\left\Vert x_{i}\right\Vert  \notag
\\
& \leq \frac{1}{2}\times \left\{ 
\begin{array}{ll}
\left( \max\limits_{1\leq i\leq n}r_{i}\right)
^{2}\sum\limits_{i=1}^{n}\left\Vert x_{i}\right\Vert ; &  \\ 
&  \\ 
\left( \sum\limits_{i=1}^{n}r_{i}^{2p}\right) ^{\frac{1}{p}}\left(
\sum\limits_{i=1}^{n}\left\Vert x_{i}\right\Vert ^{q}\right) ^{\frac{1}{q}},
& p>1,\ \frac{1}{p}+\frac{1}{q}=1; \\ 
&  \\ 
\max\limits_{1\leq i\leq n}\left\Vert x_{i}\right\Vert
\sum\limits_{i=1}^{n}r_{i}^{2}. & 
\end{array}%
\right.  \notag
\end{align}
\end{corollary}

\begin{proof}
The first part follows from Proposition \ref{p2.1} on choosing $x=x_{i},$ $%
y=e$ and applying Theorem \ref{t3.1}. The last part is obvious by H\"{o}%
lder's inequality.
\end{proof}

\begin{remark}
\label{r3.3}One would obtain the same reverse inequality (\ref{3.7}) if one
were to use Theorem \ref{t2.4}. In this case, the assumption (\ref{3.6})
should be replaced by%
\begin{equation}
\left\Vert \left\Vert x_{i}\right\Vert x_{i}-e\right\Vert \leq
r_{i}\left\Vert x_{i}\right\Vert \qquad \text{for each \qquad }i\in \left\{
1,\dots ,n\right\} .  \label{3.8}
\end{equation}
\end{remark}

On utilising the inequalities (\ref{2.5}) and (\ref{2.13.a}) one may state
the following corollary of Theorem \ref{t3.1}.

\begin{corollary}
\label{c3.4}Let $e,$ $x_{i}\in H\backslash \left\{ 0\right\} ,$ $i\in
\left\{ 1,\dots ,n\right\} $ with $\left\Vert e\right\Vert =1.$ Then we have
the inequality%
\begin{equation}
(0\leq )\sum_{i=1}^{n}\left\Vert x_{i}\right\Vert -\left\Vert
\sum_{i=1}^{n}x_{i}\right\Vert \leq \min \left\{ A,B\right\} ,  \label{3.9}
\end{equation}%
where%
\begin{equation*}
A:=2\sum_{i=1}^{n}\left\Vert x_{i}\right\Vert \left( \frac{\left\Vert
x_{i}-e\right\Vert }{\left\Vert x_{i}\right\Vert +1}\right) ^{2},
\end{equation*}%
and%
\begin{equation*}
B:=\frac{1}{2}\sum_{i=1}^{n}\frac{\left( \left\Vert x_{i}\right\Vert
+1\right) ^{2}\left\Vert x_{i}-e\right\Vert ^{2}}{\left\Vert
x_{i}\right\Vert }.
\end{equation*}
\end{corollary}

For vectors located outside the closed unit ball $\bar{B}\left( 0,1\right)
:=\left\{ z\in H|\left\Vert z\right\Vert \leq 1\right\} ,$ we may state the
following result.

\begin{corollary}
\label{c3.5}Assume that $x_{i}\notin \bar{B}\left( 0,1\right) ,$ $i\in
\left\{ 1,\dots ,n\right\} $ and $e\in H,$ $\left\Vert e\right\Vert =1.$
Then we have the inequality:%
\begin{align}
(0& \leq )\sum_{i=1}^{n}\left\Vert x_{i}\right\Vert -\left\Vert
\sum_{i=1}^{n}x_{i}\right\Vert  \label{3.10} \\
& \leq \left\{ 
\begin{array}{ll}
\dfrac{1}{2}p^{2}\sum\limits_{i=1}^{n}\left\Vert x_{i}\right\Vert
^{p-1}\left\Vert x_{i}-e\right\Vert ^{2}, & \text{if \ }p\geq 1 \\ 
&  \\ 
\dfrac{1}{2}\sum\limits_{i=1}^{n}\left\Vert x_{i}\right\Vert
^{1-p}\left\Vert x_{i}-e\right\Vert ^{2}, & \text{if \ }p<1.%
\end{array}%
\right.  \notag
\end{align}
\end{corollary}

The proof follows by Proposition \ref{p2.5.a} and Theorem \ref{t3.1}.

For complex spaces one may state the following result as well.

\begin{corollary}
\label{c3.6}Let $\left( H;\left\langle \cdot ,\cdot \right\rangle \right) $
be a complex inner product space and $\alpha _{i}\in \mathbb{C}$ with $\func{%
Re}\alpha _{i},$ $\func{Im}\alpha _{i}>0,$ $i\in \left\{ 1,\dots ,n\right\}
. $ If $x_{i},e\in H,$ $i\in \left\{ 1,\dots ,n\right\} $ with $\left\Vert
e\right\Vert =1$ and%
\begin{equation}
\left\Vert x_{i}-\frac{\func{Im}\alpha _{i}}{\func{Re}\alpha _{i}}\cdot
e\right\Vert \leq d_{i},\qquad i\in \left\{ 1,\dots ,n\right\} ,
\label{3.11}
\end{equation}%
then%
\begin{equation}
(0\leq )\sum_{i=1}^{n}\left\Vert x_{i}\right\Vert -\left\Vert
\sum_{i=1}^{n}x_{i}\right\Vert \leq \frac{1}{2}\sum_{i=1}^{n}\frac{\func{Re}%
\alpha _{i}}{\func{Im}\alpha _{i}}\cdot d_{i}^{2}.  \label{3.12}
\end{equation}
\end{corollary}

The proof follows by Theorems \ref{t2.8} and \ref{t3.1} and the details are
omitted.

Finally, by the use of Theorem \ref{t2.13}, we can state:

\begin{corollary}
\label{c3.7}If $x_{i},e\in H,$ $i\in \left\{ 1,\dots ,n\right\} $ with $%
\left\Vert e\right\Vert =1$ and $p\geq 1,$ then we have the inequalities:%
\begin{align}
(0& \leq )\sum_{i=1}^{n}\left\Vert x_{i}\right\Vert -\left\Vert
\sum_{i=1}^{n}x_{i}\right\Vert  \label{3.13} \\
& \leq \frac{1}{2}\times \left\{ 
\begin{array}{l}
\sum\limits_{i=1}^{n}\left[ \left( \left\Vert x_{i}\right\Vert +1\right)
^{2p}-\left\Vert x_{i}+e\right\Vert ^{2p}\right] ^{\frac{1}{p}}, \\ 
\\ 
\sum\limits_{i=1}^{n}\left[ \left\Vert x_{i}-e\right\Vert ^{2p}-\left\vert
\left\Vert x_{i}\right\Vert -1\right\vert ^{2p}\right] ^{\frac{1}{p}}.%
\end{array}%
\right.  \notag
\end{align}
\end{corollary}

\end{document}